\documentclass[12pt]{amsart}
\usepackage{amsmath,amsthm,epsf,amscd,amssymb,verbatim}

\newtheorem{thm}{Theorem}[section]

\theoremstyle{definition}

\theoremstyle{remark}
\newtheorem{remark}[thm]{Remark}

\numberwithin{equation}{section}

\newcommand{\bbr}{\begin{remark}}        
\newcommand{\eer}{\end{remark}}

\newcommand{\real}{\mathbb{R}}

\newcommand{\eg}{{\em e.g.}}
\newcommand{\ie}{{\em i.e.}}

\newcommand{\conf}{{\mathcal C}}

\newcommand{\dconf}{{\mathcal D}}
\newcommand{\Graph}{\Gamma}
\newcommand{\Tree}{\Upsilon}

\newcommand{\GX}{\setlength{\unitlength}{3947sp}%
\begingroup\makeatletter\ifx\SetFigFont\undefined%
\gdef\SetFigFont#1#2#3#4#5{%
  \reset@font\fontsize{#1}{#2pt}%
  \fontfamily{#3}\fontseries{#4}\fontshape{#5}%
  \selectfont}%
\fi\endgroup%
\begin{picture}(224,165)(17,657)
\thicklines
\put( 42,798){\circle{20}}
\put( 43,681){\circle{20}}
\put(216,681){\circle{20}}
\put(215,798){\circle{20}}
\put(128,739){\circle{24}}
\put( 42,798){\line( 3,-2){176.077}}
\put( 43,680){\line( 3, 2){174.692}}
\end{picture}
 }
\newcommand{\GQ}{\setlength{\unitlength}{3947sp}%
\begingroup\makeatletter\ifx\SetFigFont\undefined%
\gdef\SetFigFont#1#2#3#4#5{%
  \reset@font\fontsize{#1}{#2pt}%
  \fontfamily{#3}\fontseries{#4}\fontshape{#5}%
  \selectfont}%
\fi\endgroup%
\begin{picture}(259,148)(1,690)
\thicklines
\put(235,761){\circle{20}}
\put(136,763){\line( 1, 0){ 99}}
\put( 86,764){\circle{120}}
\put(148,764){\circle{20}}
\put( 26,763){\circle{20}}
\end{picture}
 }

\begin{document}

\title{FINDING TOPOLOGY IN A FACTORY:  CONFIGURATION SPACES}

\author{A. Abrams}
\address{Department of Mathematics, University of Georgia, 
	Athens, GA 30602, USA} 

\author{R. Ghrist}
\address{School of Mathematics and CDSNS, Georgia 
	Institute of Technology, Atlanta, GA 30332, USA}
\thanks{RG supported in part by NSF Grant \# DMS-9971629.}
\thanks{The authors wish to thank Margaret Symington for her 
	careful reading of the manuscript.}


\maketitle

It is perhaps not universally acknowledged that an outstanding 
place to find interesting topological objects is within the walls
of an automated warehouse or factory.

The examples of topological spaces constructed in this exposition 
arose simultaneously from two seemingly disparate fields:  the 
first author, in his thesis \cite{Abr00}, discovered
these spaces after working with the group of H. Landau, Z. Landau, J.
Pommersheim, and E. Zaslow \cite{TC98} on problems about multiple
random walks on graphs.
The second author \cite{GK98,G:birman} discovered these same spaces 
while collaborating with D. Koditschek in the Artificial Intelligence
Lab at the University of Michigan.  The net result makes
evident the abundance of topological 
objects within the physical world.

Topology seeks to describe, as one author puts it, the ``shape
of space'' \cite{Wee85}, with ``shape'' being interpreted as
appropriate for the context at hand. We will begin with thinking
about spaces up to homeomorphism (continuous maps with continuous
inverse), but will quickly need to abandon this class in 
favor of a looser form of equivalence: homotopy type.

Although few topological prerequisites are necessary for 
fully appreciating the examples discussed here, 
the class of spaces we consider gives an 
earthly incarnation of several intricate ideas 
from topology, such as $K(\pi,1)$ spaces (a.k.a. 
Eilenberg-MacLane spaces of type $K(\pi,1)$) and NPC
(or non-positively curved) spaces. 


\section{Configurations and Braids}

Our story begins with a classical construction:  that of a 
configuration space of points.  We consider first the 
configuration space of $N$ distinct labelled points in 
the plane $\real^2$.  That is, we consider all ordered 
$N$-tuples of points in $\real^2$ with the property that 
no two of the points 
coincide.  To each such $N$-tuple, visualized as a
collection of labeled tokens on a tabletop, we 
assign a point in the configuration space.  Two points in 
the configuration space are close if the $N$-tuples are
close as measured in $\real^{2N}$.  One can write this
mathematically as
\[ \conf^N(\real^2) := 
	\left(\real^2\times\real^2\times\cdots\times\real^2\right)
	- \Delta	,\]
where $\Delta$ denotes the ``pairwise diagonal''
\[ \Delta := \{(x_1,x_2,\ldots,x_N)\in(\real^2)^N : x_i=x_j
	{\mbox{  for some }} i\neq j\} .	\] 

This configuration space is not easy to visualize, in no small
part because the dimension of the space is $2N$. 

{\sc exercise:} Show that $\conf^2(\real^2)$ is homeomorphic 
to $\real^2\times(\real^2-(0,0))$ or equivalently 
$\real^3\times S^1$, where $S^1$ denotes the unit circle
in $\real^2$. Hint: think about placing tokens on the table one 
at a time. Does your method of proof give a simple presentation
for $\conf^3(\real^2)$? 



\begin{figure}[ht]
	{\epsfxsize=4.5in\centerline{\epsfbox{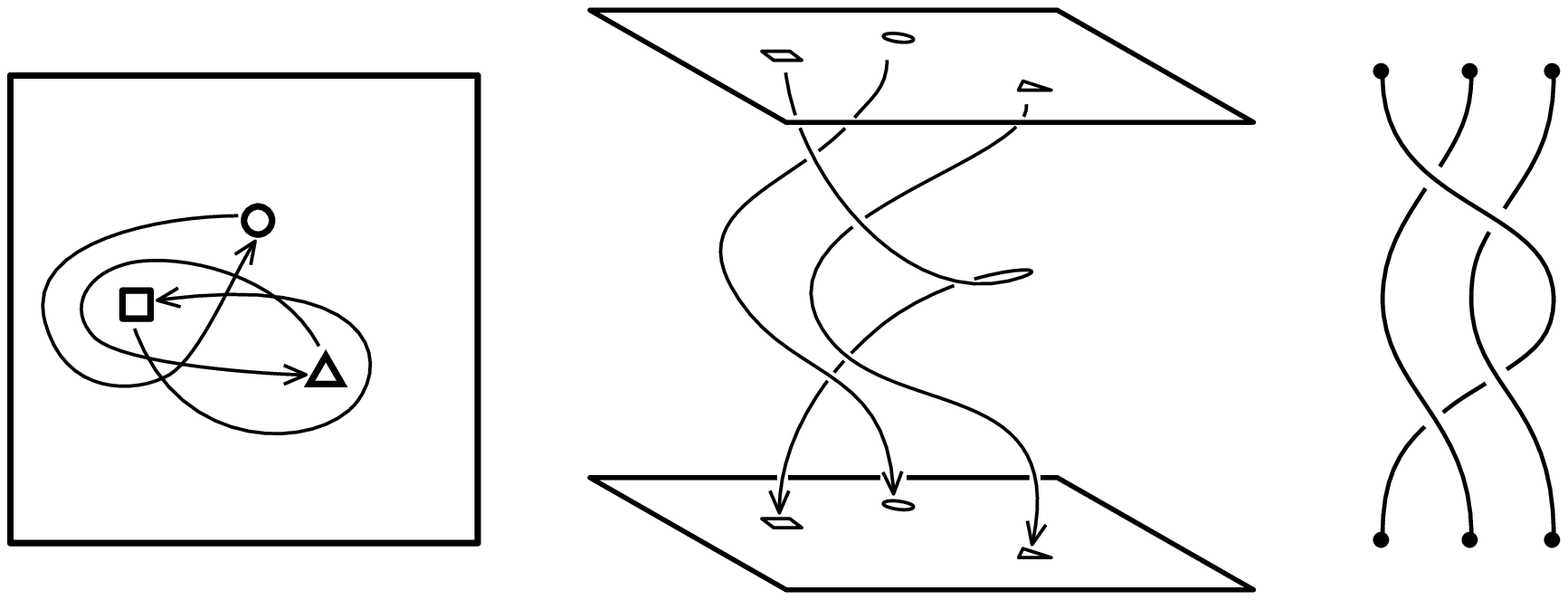}}}
	\caption{The space-time graph of a loop in $\conf^3(\real^2)$
	is a pure braid on three strands.}
	\label{fig_braid}
\end{figure}

The configuration space of points on a plane, like so many other
mathematical objects, finds multiple applications to 
physical problems. For example, any loop (closed curve) 
%
%
in $\conf^N(\real^2)$ is a {\em pure braid on $N$ 
strands}, so called because of the space-time graph of the loop
in $\conf^N(\real^2)$:  see Figure~\ref{fig_braid}. (The word 
{\em pure} indicates that each token returns to its original
location, rather than the tokens getting permuted in some other way.
Other permutations give rise to general {\em braids}, which are loops
in the space of {\em unlabeled} configurations of points in the
plane.) 
Braids are ubiquitous in the study of knots and links \cite{Bir74}, 
and thus find their way into applications ranging from quantum 
field theory to dynamical systems to DNA modeling.

\section{Robotics and Topological Motion Planning}

The space $\conf^N(\real^2)$ itself appears naturally
in robotics. 

Consider an automated factory equipped with a cadre of 
Automated Guided Vehicles (AGVs), or mobile robots, 
capable of transporting items from place to place. A 
common goal is to place several, say $N$, of these AGVs
in motion simultaneously. It is necessary to enact a 
control algorithm which perhaps moves the AGVs from 
their initial positions to their goal positions (in a 
warehousing application), or which executes a cyclic
pattern (as arises in manufacturing applications). As might be 
expected, these robots are somewhat costly and cannot tolerate
collisions without a loss of performance. Thus, modeling 
the factory floor as $\real^2$ and the AGVs as points, 
one often wishes to find paths or cycles in $\conf^N(\real^2)$  
to enact specific behaviors.\footnote{Obstacles can
easily be incorporated into these models --- a vast literature
on this subject exists \cite{latombe.book}. Executing cyclic motions
is more complex but can be at first approximated by composed
point-to-point motions. There are several assumptions required
for the simple construction we present, and various kinematic
issues (\eg, steering geometry) must in general be addressed.} 
Of course, since the robots are not truly points, and since no 
control algorithm implementation is of infinite precision, we
require that the control path reside outside of a neighborhood
of the diagonal $\Delta$ in $(\real^2)^N$.

Fortunately, this problem is not impossible to solve. 
The work of Koditschek and Rimon \cite{KR90} provides one
example of a concrete solution: they 
write out explicit vector fields on these configuration 
spaces which can be used to flow from initial to goal
positions in the presence of certain types of obstacles. 
By arranging these vector fields so that they
strongly push away from the vestiges of the diagonal $\Delta$
on the boundary of $\conf^N(\real^2)$, the control scheme is 
{\em provably} safe from collisions (as opposed to being 
{\em statistically} safe via computer simulations): no path
can ever intersect the diagonal. Furthermore, since a neighborhood
of the diagonal is repelling, the control scheme is {\em stable} 
with respect to perturbations to the system. This is quite 
important as mechanical systems have an annoying tendency to 
occasionally malfunction. Drawing the appropriate 
vector field on a configuration space yields an excellent 
method of self-correction. 

This is a clean, direct application of topological and dynamical 
ideas to a matter of great practical relevance which is currently 
used in various industrial settings.  

One reason why this problem has a nice solution is evident 
to anyone who regularly shops at a large supermarket: numerous shopping
carts trace out paths in a workspace that consists of wide
aisles. If two carts are headed toward each other, one 
needs merely steer out of the way a 
little bit (assuming the other does not move in the same direction!) 
to avoid a collision. {\em The resolution of collisions on 
$\real^2$ is a local phenomenon.}


\section{Graphs}

The robotics community, largely independently of the 
topology community, has enjoyed great success at 
identifying and manipulating configuration spaces to 
their advantage in control problems. There is, however, 
a class of simple, physically relevant scenarios
whose configuration spaces have been untapped:
the configuration spaces of points on 
a {\em graph}, or a network of edges and vertices.

Consider the situation where the AGVs must move about on a 
collection of tracks embedded in the floor, or via a 
path of electrified guide-wires from the ceiling (see 
\cite{agvhandbook} for examples). Such a restricted network is 
quite common, mainly because it is more cost-effective than a
full two-degree-of-freedom steering system for AGVs. 
In this setting, the state of the system 
at any instant of time is a point in the configuration
space of the graph $\Graph$:
\[
	\conf^N(\Graph) := 
	\left(\Graph\times\cdots\times\Graph\right)
	- \Delta	.
\]
The same principles previously mentioned still apply. To 
navigate safely on a graph, one must construct appropriate 
paths which remain strictly within $\conf^N(\Graph)$ and
are repulsed by any 
boundaries near $\Delta$. Several problems have arisen
which heretofore have prevented an analogous solution.
\begin{enumerate}
\item What do these spaces look like? 
\item How does one resolve an impending collision?
\end{enumerate}

Notice the differences between this problem and the 
problem of $\conf^N(\real^2)$. First of all, in relation 
to (1) above, $\conf^N(\Graph)$ is {\em not} a manifold:
that is, you cannot hope that every point has a neighborhood
which is locally homeomorphic to a Euclidean space. Indeed, 
if we ignore the trivial graphs which are homeomorphic to a
line segment or a circle, then the graph itself is not locally
Euclidean and products of the graph will share this feature.
Second, concerning (2) above, collisions within the 
interior of an edge are no longer locally resolvable.
Imagine that the aisles of the grocery store are only
as wide as the shopping carts, so that passing another 
person is impossible.  A store full of shoppers (using
carts) would pose a difficult
coordinated control problem.  How can carts avoid a collision
in the interior of an aisle? Clearly at least one of the 
participants 
must make a large-scale change in plans and back up to the end
of the aisle. {\em The resolution of a collision on a 
graph is a non-local phenomenon}.


\section{Examples: two robots}

Since a graph is a one-dimensional object, the configuration
space $\conf^N(\Graph)$ is $N$-dimensional. Thus configuration 
spaces of two robots are two-dimensional objects, which (at least
in simple cases) one should be able to visualize.  We therefore give
numerous examples of configuration spaces of two robots.

{\sc example 1:} $\conf^2(\setlength{\unitlength}{3947sp}%
%
\begin{picture}(221,178)(20,633)
\thicklines
\put(129,657){\circle{20}}
\put( 45,785){\circle{20}}
\put(216,785){\circle{20}}
\put(128,785){\circle{24}}
\put(129,779){\line( 0,-1){117}}
\put( 45,785){\line( 1, 0){172}}
\end{picture}
 )$

\begin{figure}[ht]
	{\epsfxsize=3.0in\centerline{\epsfbox{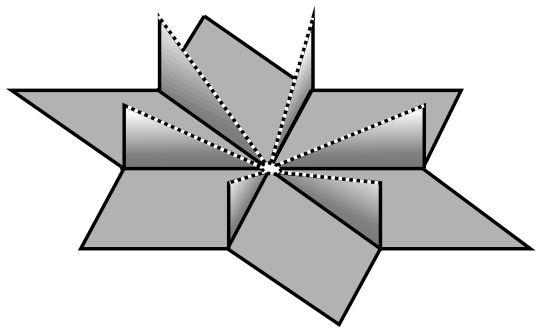}}}
	\caption{The configuration space of Example 1 embedded
	in $\real^3$. Dotted lines refer to edges which lie on the 
	diagonal $\Delta$. Note that the central vertex is deleted.}
	\label{fig_2Y}
\end{figure}

Let $\setlength{\unitlength}{3947sp}%
%
\begin{picture}(221,178)(20,633)
\thicklines
\put(129,657){\circle{20}}
\put( 45,785){\circle{20}}
\put(216,785){\circle{20}}
\put(128,785){\circle{24}}
\put(129,779){\line( 0,-1){117}}
\put( 45,785){\line( 1, 0){172}}
\end{picture}
 $ denote the graph of three edges attached at a 
central vertex. The cellular structure of $\conf^2(\setlength{\unitlength}{3947sp}%
%
\begin{picture}(221,178)(20,633)
\thicklines
\put(129,657){\circle{20}}
\put( 45,785){\circle{20}}
\put(216,785){\circle{20}}
\put(128,785){\circle{24}}
\put(129,779){\line( 0,-1){117}}
\put( 45,785){\line( 1, 0){172}}
\end{picture}
 )$ is 
simple to procure. There are six (two robots times 
three-choose-two edges) squares or ``2-cells'', each with one 
corner punctured, corresponding to
those configurations where the two robots are on distinct
edges of $\setlength{\unitlength}{3947sp}%
%
\begin{picture}(221,178)(20,633)
\thicklines
\put(129,657){\circle{20}}
\put( 45,785){\circle{20}}
\put(216,785){\circle{20}}
\put(128,785){\circle{24}}
\put(129,779){\line( 0,-1){117}}
\put( 45,785){\line( 1, 0){172}}
\end{picture}
 $. Since there are three edges in \setlength{\unitlength}{3947sp}%
%
\begin{picture}(221,178)(20,633)
\thicklines
\put(129,657){\circle{20}}
\put( 45,785){\circle{20}}
\put(216,785){\circle{20}}
\put(128,785){\circle{24}}
\put(129,779){\line( 0,-1){117}}
\put( 45,785){\line( 1, 0){172}}
\end{picture}
 , those 
remaining configurations in which both robots are on the same 
edge yield three square cells, each of which is divided by
the diagonal $\Delta$ into a pair of triangular cells. 
Thus there are six triangular 2-cells
corresponding to those configurations where both robots are on 
the same edge, but distinct. By enumerating the behaviors of
each of these 2-cells, one can make the identifications to 
arrive at the space given in Figure~\ref{fig_2Y} (or, one can 
use a slightly more sophisticated argument as in 
\cite{GK98,GK00}).

{\sc exercise:} Place two coins on a piece of paper with a large
\setlength{\unitlength}{3947sp}%
%
\begin{picture}(221,178)(20,633)
\thicklines
\put(129,657){\circle{20}}
\put( 45,785){\circle{20}}
\put(216,785){\circle{20}}
\put(128,785){\circle{24}}
\put(129,779){\line( 0,-1){117}}
\put( 45,785){\line( 1, 0){172}}
\end{picture}
  drawn on it so that both coins are on the same edge of
the graph. Using two fingers, execute a path which exchanges 
the coins' positions without collisions. Draw this motion as
a path on Figure~\ref{fig_2Y}. 

{\sc example 2:} $\conf^2(\setlength{\unitlength}{3947sp}%
\begingroup\makeatletter\ifx\SetFigFont\undefined%
\gdef\SetFigFont#1#2#3#4#5{%
  \reset@font\fontsize{#1}{#2pt}%
  \fontfamily{#3}\fontseries{#4}\fontshape{#5}%
  \selectfont}%
\fi\endgroup%
\begin{picture}(259,148)(1,690)
\thicklines
\put(235,761){\circle{20}}
\put(136,763){\line( 1, 0){ 99}}
\put( 86,764){\circle{120}}
\put(148,764){\circle{20}}
\put( 26,763){\circle{20}}
\end{picture}
 )$

\begin{figure}[ht]
	{\epsfxsize=4.5in\centerline{\epsfbox{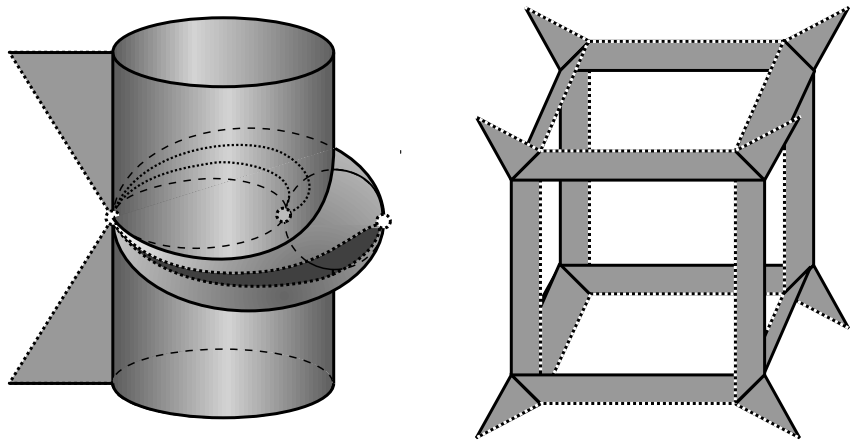 }}}
	\caption{The configuration spaces of Examples 2 [left]
	and 3 [right]. Dotted lines refer to 
	vestiges of the diagonal $\Delta$.}
	\label{fig_2QX}
\end{figure}

Let $\setlength{\unitlength}{3947sp}%
\begingroup\makeatletter\ifx\SetFigFont\undefined%
\gdef\SetFigFont#1#2#3#4#5{%
  \reset@font\fontsize{#1}{#2pt}%
  \fontfamily{#3}\fontseries{#4}\fontshape{#5}%
  \selectfont}%
\fi\endgroup%
\begin{picture}(259,148)(1,690)
\thicklines
\put(235,761){\circle{20}}
\put(136,763){\line( 1, 0){ 99}}
\put( 86,764){\circle{120}}
\put(148,764){\circle{20}}
\put( 26,763){\circle{20}}
\end{picture}
 $ denote the graph with three edges obtained from 
$\setlength{\unitlength}{3947sp}%
%
\begin{picture}(221,178)(20,633)
\thicklines
\put(129,657){\circle{20}}
\put( 45,785){\circle{20}}
\put(216,785){\circle{20}}
\put(128,785){\circle{24}}
\put(129,779){\line( 0,-1){117}}
\put( 45,785){\line( 1, 0){172}}
\end{picture}
 $ by gluing two boundary vertices together. 
One method of constructing 
$\conf^2(\setlength{\unitlength}{3947sp}%
\begingroup\makeatletter\ifx\SetFigFont\undefined%
\gdef\SetFigFont#1#2#3#4#5{%
  \reset@font\fontsize{#1}{#2pt}%
  \fontfamily{#3}\fontseries{#4}\fontshape{#5}%
  \selectfont}%
\fi\endgroup%
\begin{picture}(259,148)(1,690)
\thicklines
\put(235,761){\circle{20}}
\put(136,763){\line( 1, 0){ 99}}
\put( 86,764){\circle{120}}
\put(148,764){\circle{20}}
\put( 26,763){\circle{20}}
\end{picture}
 )$ would be to first remove the
configurations where both robots are on the vertices to be
glued. Then identify those portions of the 
boundary of $\conf^2(\setlength{\unitlength}{3947sp}%
%
\begin{picture}(221,178)(20,633)
\thicklines
\put(129,657){\circle{20}}
\put( 45,785){\circle{20}}
\put(216,785){\circle{20}}
\put(128,785){\circle{24}}
\put(129,779){\line( 0,-1){117}}
\put( 45,785){\line( 1, 0){172}}
\end{picture}
 )$ which have a robot at the 
vertices to be glued in $\setlength{\unitlength}{3947sp}%
\begingroup\makeatletter\ifx\SetFigFont\undefined%
\gdef\SetFigFont#1#2#3#4#5{%
  \reset@font\fontsize{#1}{#2pt}%
  \fontfamily{#3}\fontseries{#4}\fontshape{#5}%
  \selectfont}%
\fi\endgroup%
\begin{picture}(259,148)(1,690)
\thicklines
\put(235,761){\circle{20}}
\put(136,763){\line( 1, 0){ 99}}
\put( 86,764){\circle{120}}
\put(148,764){\circle{20}}
\put( 26,763){\circle{20}}
\end{picture}
 $, and glue these
portions of $\conf^2(\setlength{\unitlength}{3947sp}%
%
\begin{picture}(221,178)(20,633)
\thicklines
\put(129,657){\circle{20}}
\put( 45,785){\circle{20}}
\put(216,785){\circle{20}}
\put(128,785){\circle{24}}
\put(129,779){\line( 0,-1){117}}
\put( 45,785){\line( 1, 0){172}}
\end{picture}
 )$ together. The result, although a 
very simple configuration space, is already somewhat
difficult to visualize: we illustrate the space, embedded in 
$\real^3$, in Figure~\ref{fig_2QX}[left]. There are 
three ``punctures'' at which both robots collide at
one of the three vertices. The six dotted edges are
the images of the diagonal curves from Figure~\ref{fig_2Y}.

{\sc example 3:} $\conf^2(\setlength{\unitlength}{3947sp}%
\begingroup\makeatletter\ifx\SetFigFont\undefined%
\gdef\SetFigFont#1#2#3#4#5{%
  \reset@font\fontsize{#1}{#2pt}%
  \fontfamily{#3}\fontseries{#4}\fontshape{#5}%
  \selectfont}%
\fi\endgroup%
\begin{picture}(224,165)(17,657)
\thicklines
\put( 42,798){\circle{20}}
\put( 43,681){\circle{20}}
\put(216,681){\circle{20}}
\put(215,798){\circle{20}}
\put(128,739){\circle{24}}
\put( 42,798){\line( 3,-2){176.077}}
\put( 43,680){\line( 3, 2){174.692}}
\end{picture}
 )$

Increasing the incidence number of the central vertex 
complicates the configuration space. Consider 
$\setlength{\unitlength}{3947sp}%
\begingroup\makeatletter\ifx\SetFigFont\undefined%
\gdef\SetFigFont#1#2#3#4#5{%
  \reset@font\fontsize{#1}{#2pt}%
  \fontfamily{#3}\fontseries{#4}\fontshape{#5}%
  \selectfont}%
\fi\endgroup%
\begin{picture}(224,165)(17,657)
\thicklines
\put( 42,798){\circle{20}}
\put( 43,681){\circle{20}}
\put(216,681){\circle{20}}
\put(215,798){\circle{20}}
\put(128,739){\circle{24}}
\put( 42,798){\line( 3,-2){176.077}}
\put( 43,680){\line( 3, 2){174.692}}
\end{picture}
 $ a radial tree of four edges emanating from a 
central vertex. The visualization of $\conf^2(\setlength{\unitlength}{3947sp}%
\begingroup\makeatletter\ifx\SetFigFont\undefined%
\gdef\SetFigFont#1#2#3#4#5{%
  \reset@font\fontsize{#1}{#2pt}%
  \fontfamily{#3}\fontseries{#4}\fontshape{#5}%
  \selectfont}%
\fi\endgroup%
\begin{picture}(224,165)(17,657)
\thicklines
\put( 42,798){\circle{20}}
\put( 43,681){\circle{20}}
\put(216,681){\circle{20}}
\put(215,798){\circle{20}}
\put(128,739){\circle{24}}
\put( 42,798){\line( 3,-2){176.077}}
\put( 43,680){\line( 3, 2){174.692}}
\end{picture}
 )$ is 
a bit more involved and requires some work to obtain. 
For the purpose of stimulating curiosity, we include this
configuration space as Figure~\ref{fig_2QX}[right]. 


\section{Simplification: discretization}

To visualize these configuration spaces, it is clear that 
some simplification into a more manageable form is
necessary. There are two principal methods we use:
the first is a way of removing the ``unsafe'' cells of the 
space near the diagonal $\Delta$. The second, deformation 
retraction, is a more drastic
crushing of the space down to a lower dimensional 
``skeleton,'' as discussed in the subsequent section.

Any graph $\Graph$ comes equipped with a cellular structure:
0-cells (vertices) and 1-cells (edges). The $N$-fold cross
product of $\Graph$ with itself inherits a cell structure, 
each cell being a product of $N$ (not necessarily distinct) cells
in $\Graph$: cf. Example 1. However, the configuration space
does not quite have a natural cell structure, since the diagonal
$\Delta$ slices through all product cells with repeated factors.
Notice, however, that in several of the previous examples, these 
partial cells dangle ``inessentially'' and could be collapsed onto
a more ``essential'' skeleton of the configuration space.

Such an operation can be interpreted as follows \cite{Abr00}.
Consider the {\em discretized configuration space} of $\Graph$, 
denoted $\dconf^N(\Graph)$, defined as 
$\left(\Graph\times\cdots\times\Graph\right)-\tilde\Delta$, 
where $\tilde\Delta$ denotes the set of all product cells 
in $\Graph\times\cdots\times\Graph$ whose closures intersect
the diagonal $\Delta$. An equivalent description of 
$\dconf^N(\Graph)$
is that it is the set of configurations for which, given any two 
robots on $\Graph$ and any path in $\Graph$ connecting them, the 
path contains at least one entire edge. 
Thus, instead of restricting robots to be at least some intrinsic
distance $\epsilon$ apart (i.e., removing an $\epsilon$ neighborhood
of $\Delta$), one now restricts robots on $\Graph$ to be ``at least
one full edge apart.'' This is a natural kind of configuration
space in the context of random walks on graphs \cite{TC98}. Note 
that $\dconf^N(\Graph)$ is a subcomplex of $\conf^N(\Graph)$
(it does not contain ``partial cells'' which arise when cutting
along the diagonal), and is, in fact, the largest subcomplex of
$\conf^N(\Graph)$ which does not intersect $\Delta$. 

With this natural cell structure, one can think of the vertices 
(0-cells) of $\dconf^N(\Graph)$ as ``discretized'' configurations 
--- arrangements of 
labeled tokens at the vertices of the graph. The edges of
$\dconf^N(\Graph)$, or 1-cells, tell you which discrete 
configurations can be connected by moving one token along an
edge of $\Graph$. Each 2-cell in $\dconf^N(\Graph)$ 
represents two independent (or ``commuting'') edges:
one can move a pair of tokens independently along disjoint
edges. A $k$-cell in $\dconf^N(\Graph)$ likewise represents 
the ability to move $k$ tokens along $k$ disjoint edges
in $\Graph$.

Returning to Figure~\ref{fig_2Y}, the discretization of $\conf^2(\setlength{\unitlength}{3947sp}%
%
\begin{picture}(221,178)(20,633)
\thicklines
\put(129,657){\circle{20}}
\put( 45,785){\circle{20}}
\put(216,785){\circle{20}}
\put(128,785){\circle{24}}
\put(129,779){\line( 0,-1){117}}
\put( 45,785){\line( 1, 0){172}}
\end{picture}
 )$ 
removes much of the space. For example, the triangular two-dimensional
cells represent configurations where both robots are on the interior
of the same edge. Since these are not ``one full edge apart,'' these
cells are deleted. The same is true of all the other two-dimensional 
cells which represent robots in the interior of separate edges. Which 
configurations of two robots on $\setlength{\unitlength}{3947sp}%
%
\begin{picture}(221,178)(20,633)
\thicklines
\put(129,657){\circle{20}}
\put( 45,785){\circle{20}}
\put(216,785){\circle{20}}
\put(128,785){\circle{24}}
\put(129,779){\line( 0,-1){117}}
\put( 45,785){\line( 1, 0){172}}
\end{picture}
 $ are separated by a full edge?

\begin{figure}[ht]
	{\epsfxsize=4.75in\centerline{\epsfbox{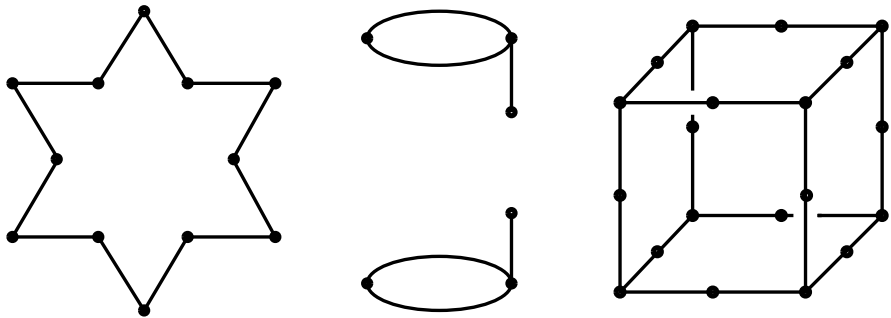 }}}
	\caption{The discretizations of the configuration spaces
	in Examples 1-3 (left to right). }
	\label{fig_discs}
\end{figure}

{\sc exercise:} Show that discretizing the configuration spaces 
of Examples 1 through 3 yields the configuration spaces of 
Figure~\ref{fig_discs}. How good of an ``approximation'' are
these spaces?

One could compute the discretization of Example 1 in a less direct 
manner which generalizes to some lovely examples to follow.
Recall that the discretized
configuration spaces inherit a cell structure from $\Graph$: all 
the product cells of $\Graph^N$ not entirely in $\conf^N(\Graph)$ are
removed by the discretization. Thus, in the case of Example 1,
simple counting reveals that the space $\dconf^2(\setlength{\unitlength}{3947sp}%
%
\begin{picture}(221,178)(20,633)
\thicklines
\put(129,657){\circle{20}}
\put( 45,785){\circle{20}}
\put(216,785){\circle{20}}
\put(128,785){\circle{24}}
\put(129,779){\line( 0,-1){117}}
\put( 45,785){\line( 1, 0){172}}
\end{picture}
 )$ possesses 
twelve 0-cells (where both robots are at distinct vertices of $\setlength{\unitlength}{3947sp}%
%
\begin{picture}(221,178)(20,633)
\thicklines
\put(129,657){\circle{20}}
\put( 45,785){\circle{20}}
\put(216,785){\circle{20}}
\put(128,785){\circle{24}}
\put(129,779){\line( 0,-1){117}}
\put( 45,785){\line( 1, 0){172}}
\end{picture}
 $), 
twelve 1-cells (where one robot is at a vertex and the other is on an 
edge whose closure does not contain said vertex), and zero 2-cells (since
every pair of edges intersect along their boundaries). With a 
little thought, one can see that $\dconf^2(\setlength{\unitlength}{3947sp}%
%
\begin{picture}(221,178)(20,633)
\thicklines
\put(129,657){\circle{20}}
\put( 45,785){\circle{20}}
\put(216,785){\circle{20}}
\put(128,785){\circle{24}}
\put(129,779){\line( 0,-1){117}}
\put( 45,785){\line( 1, 0){172}}
\end{picture}
 )$ is a connected 
manifold: each zero-cell connects to exactly two 1-cells, and all
of the 1-cells are joined end-to-end cyclically. Thus,
$\dconf^2(\setlength{\unitlength}{3947sp}%
%
\begin{picture}(221,178)(20,633)
\thicklines
\put(129,657){\circle{20}}
\put( 45,785){\circle{20}}
\put(216,785){\circle{20}}
\put(128,785){\circle{24}}
\put(129,779){\line( 0,-1){117}}
\put( 45,785){\line( 1, 0){172}}
\end{picture}
 )$ is a topological circle, precisely as
obtained by deleting all the near-diagonal cells from 
$\conf^2(\setlength{\unitlength}{3947sp}%
%
\begin{picture}(221,178)(20,633)
\thicklines
\put(129,657){\circle{20}}
\put( 45,785){\circle{20}}
\put(216,785){\circle{20}}
\put(128,785){\circle{24}}
\put(129,779){\line( 0,-1){117}}
\put( 45,785){\line( 1, 0){172}}
\end{picture}
 )$ in Figure~\ref{fig_2Y}. The discretization operation 
yields a subcomplex of $\conf^2(\setlength{\unitlength}{3947sp}%
%
\begin{picture}(221,178)(20,633)
\thicklines
\put(129,657){\circle{20}}
\put( 45,785){\circle{20}}
\put(216,785){\circle{20}}
\put(128,785){\circle{24}}
\put(129,779){\line( 0,-1){117}}
\put( 45,785){\line( 1, 0){172}}
\end{picture}
 )$ which appears to contain all the 
``essential'' topology (more specifically, the spaces $\conf^N(\setlength{\unitlength}{3947sp}%
%
\begin{picture}(221,178)(20,633)
\thicklines
\put(129,657){\circle{20}}
\put( 45,785){\circle{20}}
\put(216,785){\circle{20}}
\put(128,785){\circle{24}}
\put(129,779){\line( 0,-1){117}}
\put( 45,785){\line( 1, 0){172}}
\end{picture}
 )$ 
and $\dconf^N(\setlength{\unitlength}{3947sp}%
%
\begin{picture}(221,178)(20,633)
\thicklines
\put(129,657){\circle{20}}
\put( 45,785){\circle{20}}
\put(216,785){\circle{20}}
\put(128,785){\circle{24}}
\put(129,779){\line( 0,-1){117}}
\put( 45,785){\line( 1, 0){172}}
\end{picture}
 )$ are of the same {\em homotopy type} --- see the 
next section for definitions); however, this is 
certainly not the case for the discretization of $\conf^2(\setlength{\unitlength}{3947sp}%
\begingroup\makeatletter\ifx\SetFigFont\undefined%
\gdef\SetFigFont#1#2#3#4#5{%
  \reset@font\fontsize{#1}{#2pt}%
  \fontfamily{#3}\fontseries{#4}\fontshape{#5}%
  \selectfont}%
\fi\endgroup%
\begin{picture}(259,148)(1,690)
\thicklines
\put(235,761){\circle{20}}
\put(136,763){\line( 1, 0){ 99}}
\put( 86,764){\circle{120}}
\put(148,764){\circle{20}}
\put( 26,763){\circle{20}}
\end{picture}
 )$, 
which becomes disconnected! In the next section, we will state
the criteria under which discretization is topologically faithful.

The counting arguments used above often can determine the discretized 
configuration space, even when the full configuration space is
unknown. The following are some surprising examples of interesting 
spaces arising as the discretized configuration space of non-planar
graphs \cite{Abr00}. 
%
%
%
%

\begin{figure}[ht]
	{\epsfxsize=3.25in\centerline{\epsfbox{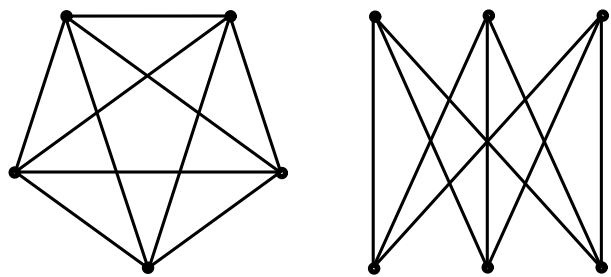 }}}
	\caption{The non-planar graphs $K_{5}$ [left] and 
	$K_{3,3}$ [right]. This notation comes from graph
	theory, where these are fundamental examples of 
	non-planar graphs.}
	\label{fig_nonplanar}
\end{figure}

{\sc example 4:} $\dconf^2(K_{5})$

Consider the complete bipartite graph $K_{5}$ pictured in 
Figure~\ref{fig_nonplanar}[left]. 
The discretized configuration space of two robots on this graph 
is a two-dimensional complex. A simple counting argument reveals 
the cell-structure:
\begin{figure}[ht]
	{\epsfxsize=5.4in\epsfbox{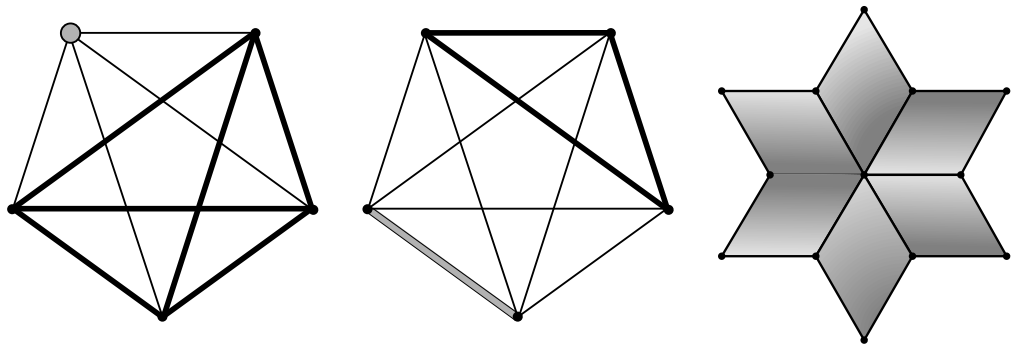}}
	\caption{[left] For every vertex in the space 
	$K_{5}$ there are six disjoint edges.
	Likewise [middle] for each edge there are three
	totally disjoint edges. In $\dconf^2(K_5)$, these cells 
	fit together to form a locally Euclidean 
	two-dimensional complex [right].}
	\label{fig_cells}
\end{figure}
\begin{enumerate} 
\item[{\sc 0-cells}] Each 0-cell corresponds to a configuration in which 
	the two robots are at distinct vertices. Since $K_{5}$
	has five vertices, there are exactly $(5)(5-1)=20$ such 
	0-cells.  (Remember, there is no vertex where two edges cross
	in the picture; there are only vertices at the corners
	of the pentagon.)
\item[{\sc 1-cells}] Each 1-cell corresponds to a configuration in which 
	one robot is at a vertex and the other is on an edge 
	whose endpoints do not include the vertex already
	occupied. From the diagram of $K_{5}$ one counts
	that there are $(2)(5)(6)=60$ such 1-cells, as in
	Figure~\ref{fig_cells}[left]. The factor of 
	two comes from the fact that we label the two  
	robots on $K_{5}$.
\item[{\sc 2-cells}] Each 2-cell corresponds to a configuration in which
	the two robots occupy edges whose closures are disjoint.
	Again, from the diagram (and Figure~\ref{fig_cells}[middle])
	one counts that there are $(10)(3)=30$ 
	such 2-cells in the complex. 
\end{enumerate}
One then demonstrates that each edge borders a pair of 2-cells 
preserving an orientation and that each vertex is incident to 
six edges, as in Figure~\ref{fig_cells}[right].  Also, the space 
$\dconf^2(K_{5})$ is
connected: it is not hard to see that you can move from any
configuration to any other. Thus $\dconf^2(K_{5})$ is a connected 
orientable surface, and the Classification Theorem for surfaces 
implies that the space is determined uniquely up to homeomorphism
by the {\em Euler characteristic},
\begin{equation}
	\chi(\dconf^2(K_{5})) := \#{\mbox{faces}}
				-\#{\mbox{edges}}
				+\#{\mbox{vertices}} = 30-60+20 = -10 .
\end{equation}
Thus, $\dconf^2(K_{5})$ is a closed orientable surface of 
genus $g:=1-\frac{1}{2}\chi = 6$. It is not at all obvious that 
the motion of two robots on this graph should produce a 
genus six surface. Obtaining a manifold is surprising enough,
but a manifold with genus larger than one really goes against
the notion that all of the interesting topology in these
spaces is ``localized'' in configurations about a vertex. 

\begin{figure}[ht]
	{\epsfxsize=5.0in\centerline{\epsfbox{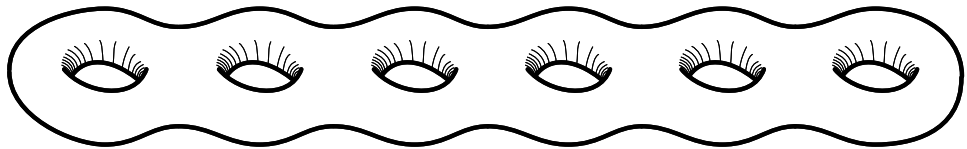 }}}
	\caption{The space $\dconf^2(K_5)$ is homeomorphic to
	a closed orientable surface of genus six.}
	\label{fig_genus6}
\end{figure}

{\sc example 5:} $\dconf^2(K_{3,3})$

A near-identical analysis on the graph $K_{3,3}$ of 
Figure~\ref{fig_nonplanar}[right] reveals that $\dconf^2(K_{3,3})$ 
is also a connected closed orientable surface. 
The natural cell structure on $\dconf^2(K_{3,3})$ possesses 
exactly 36 faces, 72 edges, and 30 vertices. Thus, 
\begin{equation}
\chi(\dconf^2(K_5)) = 36-72+30=-6,
\end{equation} 
and we conclude that the genus of this surface is four.

In order to understand discretizations of higher-dimensional 
configuration spaces, one can employ an appropriate version of 
duality \cite{Abr00}: instead of tracking how distinct robots on 
$\Graph$ move about, one considers configurations of ``holes'' ---
regions on $\Graph$ which have no robots on them. As the holes
have no natural labeling, the duality applies directly to the 
{\em unlabeled} configuration spaces.  It nevertheless implies
a relationship between $\dconf^N(\Graph)$ and $\dconf^{V-N}(\Graph)$,
where $V$ is the number of vertices
of $\Graph$. Using this duality, one obtains the following 
examples.

{\sc example 6:} The space $\dconf^3(K_{5})$ is homeomorphic
to a connected closed orientable surface of genus 16.

{\sc example 7:} Likewise, $\dconf^4(K_{3,3})$ is
homeomorphic to a connected closed orientable surface of genus 37.


\section{Simplification: deformation}

We are naturally led to the question of how good an approximation 
the discretized configuration space is. Heuristically, the spaces
$\conf^N$ and $\dconf^N$ should be ``similar,'' since the latter
is a subset of the former obtained by collapsing out those 
cells which border the remains of the diagonal $\Delta$. That the
discretization is not always faithful is evidenced by the
fact that $\conf^2(\setlength{\unitlength}{3947sp}%
\begingroup\makeatletter\ifx\SetFigFont\undefined%
\gdef\SetFigFont#1#2#3#4#5{%
  \reset@font\fontsize{#1}{#2pt}%
  \fontfamily{#3}\fontseries{#4}\fontshape{#5}%
  \selectfont}%
\fi\endgroup%
\begin{picture}(259,148)(1,690)
\thicklines
\put(235,761){\circle{20}}
\put(136,763){\line( 1, 0){ 99}}
\put( 86,764){\circle{120}}
\put(148,764){\circle{20}}
\put( 26,763){\circle{20}}
\end{picture}
 )$ is connected while $\dconf^2(\setlength{\unitlength}{3947sp}%
\begingroup\makeatletter\ifx\SetFigFont\undefined%
\gdef\SetFigFont#1#2#3#4#5{%
  \reset@font\fontsize{#1}{#2pt}%
  \fontfamily{#3}\fontseries{#4}\fontshape{#5}%
  \selectfont}%
\fi\endgroup%
\begin{picture}(259,148)(1,690)
\thicklines
\put(235,761){\circle{20}}
\put(136,763){\line( 1, 0){ 99}}
\put( 86,764){\circle{120}}
\put(148,764){\circle{20}}
\put( 26,763){\circle{20}}
\end{picture}
 )$ is
not. Indeed, by definition, $\dconf^N(\Graph)$ is the empty set
whenever $N$ is greater than the number of vertices of $\Graph$: 
\ie, when the discretization of $\Graph$ is too ``coarse.''

The notion of ``sameness'' appropriate here is that of {\em
deformation retract}. A subspace $A$ of a space $X$ is a (strong) 
deformation retract of $X$ if there exists a continuous family of 
continuous maps $f_t:X\to X$ such that $f_0$ is the identity map on $X$ 
and $f_1$ is a map which sends $X$ onto $A$, and such that $f_t$ fixes 
$A$ pointwise for all $t$. The image of each $f_t$ can be seen as
frames in a movie which exhibits a continuous shrinking of 
$X$ onto $A$. Deformation retractions
are an excellent way to simplify a space without changing any
essential topological properties.\footnote{Except perhaps dimension, 
but for the applications at hand, reducing the dimension of the
configuration space is a boon.} The most fundamental notion of 
topological equivalence, {\em homotopy type}, can be defined
in terms of deformation retracts: two spaces $X$ and $Y$ are
of the same homotopy type if and only if they are both deformation 
retracts of a ``larger'' space $Z$. 

The key result is that $\dconf^N(\Graph)$ is a deformation retract 
of $\conf^N(\Graph)$ as long as the discretization of the graph
$\Graph$ is not too coarse. More specifically, 

{\bf Theorem 1:} [Abrams \cite{Abr00}] {\em
For any $N>1$ and any graph $\Graph$ with at least $N$ vertices, 
$\conf^N(\Graph)$ deformation retracts to $\dconf^N(\Graph)$ 
if and only if
\begin{enumerate}
\item Each path between distinct vertices of valence not equal to two 
	passes through at least $N-1$ edges; and
\item Each loop from a vertex to itself which cannot be 
	shrunk to a point in $\Graph$ passes through at least
	$N+1$ edges.
\end{enumerate}
}
(The valence of a vertex is the number of incident edges.)

It thus follows that the spaces $\conf^2(K_{5})$ and 
$\conf^2(K_{3,3})$ deformation retract to the discretized 
configuration spaces computed in Examples 4 and 5. This is
extremely useful information: trying to compute the 
configuration space $\conf^2(K_5)$ directly would appear
hopelessly complex. Note 
that the second condition fails for the discretization 
of $\conf^2(\setlength{\unitlength}{3947sp}%
\begingroup\makeatletter\ifx\SetFigFont\undefined%
\gdef\SetFigFont#1#2#3#4#5{%
  \reset@font\fontsize{#1}{#2pt}%
  \fontfamily{#3}\fontseries{#4}\fontshape{#5}%
  \selectfont}%
\fi\endgroup%
\begin{picture}(259,148)(1,690)
\thicklines
\put(235,761){\circle{20}}
\put(136,763){\line( 1, 0){ 99}}
\put( 86,764){\circle{120}}
\put(148,764){\circle{20}}
\put( 26,763){\circle{20}}
\end{picture}
 )$, but adding one more vertex would yield a 
faithful discretized configuration space. The 
discretizations of Examples 6 and 7 are not fine enough to 
give an equivalence. 

The dimension of the smallest subcomplex to which a configuration
space deformation retracts is an important quantity in 
practice, since a large dimension greatly increases the 
complexity of the computational work needed to control the
system. The following theorem reveals that the essential 
dimension of the configuration space is governed by 
properties of the graph, independent of the number of 
robots on the graph.

{\bf Theorem 2:} [Ghrist \cite{G:birman}] {\em
Given a graph $\Graph$ having $V$ vertices of
valence greater than two, the space $\conf^N(\Graph)$ 
deformation retracts to a subcomplex of dimension at most $V$.
}

{\sc example 8:} $\conf^N(\Tree_k)$

Consider the radial $k$-prong tree $\Tree_k$ having $k>2$ edges 
and $k+1$ vertices, all edges being attached at a single central 
vertex. For example, $\Tree_3=\setlength{\unitlength}{3947sp}%
%
\begin{picture}(221,178)(20,633)
\thicklines
\put(129,657){\circle{20}}
\put( 45,785){\circle{20}}
\put(216,785){\circle{20}}
\put(128,785){\circle{24}}
\put(129,779){\line( 0,-1){117}}
\put( 45,785){\line( 1, 0){172}}
\end{picture}
 $, and $\Tree_4=\setlength{\unitlength}{3947sp}%
\begingroup\makeatletter\ifx\SetFigFont\undefined%
\gdef\SetFigFont#1#2#3#4#5{%
  \reset@font\fontsize{#1}{#2pt}%
  \fontfamily{#3}\fontseries{#4}\fontshape{#5}%
  \selectfont}%
\fi\endgroup%
\begin{picture}(224,165)(17,657)
\thicklines
\put( 42,798){\circle{20}}
\put( 43,681){\circle{20}}
\put(216,681){\circle{20}}
\put(215,798){\circle{20}}
\put(128,739){\circle{24}}
\put( 42,798){\line( 3,-2){176.077}}
\put( 43,680){\line( 3, 2){174.692}}
\end{picture}
 $. 
From the above theorem, we see that $\conf^N(\Tree_k)$
must deformation retract to a 1-dimensional subcomplex --- that is, 
a graph. Since the essential topological features of a graph
are determined by its Euler characteristic, one need merely
compute the number of vertices and edges to classify these spaces.
Using a double-induction argument on $N$ and $k$ \cite{G:birman}, 
one derives a two-variable recursion relation for the Euler 
characteristic. By solving this equation, one can prove that 
$\conf^N(\Tree_k)$ has the homotopy type of a graph which is
a ``bouquet'' of $P$ distinct loops joined together like petals
on a daisy, where
\begin{equation}
P = 1 + (Nk-2N-k+1)\frac{(N+k-2)!}{(k-1)!} .
\end{equation}
Note, for example, that $\conf^2(\Tree_3)$ has exactly one generating
loop, as Figures~\ref{fig_2Y} and \ref{fig_discs}[left] confirm.
Figure~\ref{fig_discs}[right] provides another confirmation of 
this equation for $k=4$, as the reader should verify.
The factorial growth of $P$ in $N$ is due to the 
fact that we label the $N$ robots on $\Tree_k$. If one considers
the unlabeled configuration spaces, then
the second term in the expression above is reduced 
by a factor of $N!$.

It is worth emphasizing that while the control problem of 
robots on a graph is rather intuitive for two robots, it quickly
builds in complexity. Since the dimension alone makes the 
configuration spaces in general nearly impossible to visualize,
Theorem 2 is quite helpful --- the ``essential'' dimension of
the configuration space is independent of the number of robots
on the graph. For the graph $\Tree_k$, Theorem 2 implies that there
is a one-dimensional {\em roadmap} which gives a  
perfect representation of the configuration space: no topological
data is lost. Since the proof of Theorem 2 is constructive, 
one can use standard algorithms for determining shortest paths
on a graph in order to develop efficient path planning for 
multiple robots on $\Tree_k$ via the roadmap. 


\section{Conclusions}

Applications of configuration spaces to robotics are by no means 
novel: ideas by mathematicians, computer scientists, 
and engineers have been steadily developing since the 
1960's (see, e.g., \cite{Can88,Got88} for an introduction).
In addition, various kinds of configuration spaces 
arise in topology and physics rather often, as in the study of braids
\cite{Bir74}, abstract linkages \cite{Wal85,KM96}, or
invariants of manifolds \cite{BottTaub}.

In the particular application involving multiple independent robots, 
the global aspect of the control problem is the principal difficulty 
when the robots are constrained to a network. It is this global 
nature which hints at the efficacy of a topological viewpoint. 
Indeed, the determination and simplification of the configuration
spaces of graphs can be used to construct practical control schemes
(see \cite{GK98,GK00} for some simple examples). 

Reversing the perspective, it is remarkable that this class
of topologically rich configuration spaces were virtually 
untouched until motivated by problems from other fields. There
are many deeper properties of these spaces which one can 
prove \cite{G:birman,Abr00}. We list a few below:
\begin{enumerate}
\item
	For any graph $\Graph$, $\conf^N(\Graph)$ is an 
	Eilenberg-MacLane space of type $K(\pi,1)$. That
	is, the image of any continuous map from a $k$-dimensional
	sphere $S^k$, $k>1$, into $\conf^N(\Graph)$ can be 
	shrunk to a point in $\conf^N(\Graph)$. Such a space 
	is sometimes descriptively called ``aspherical.''
\item 
	The discretized configuration space has a natural
	structure of a cube-complex, since all the cells are
	products of intervals. This allows one to use recent
	fast algorithms from computational homology 
	to determine homology groups and 
	generators in practical settings.
\item
	From the cube complex structure of $\dconf^N(\Graph)$, 
	one can show that $\conf^N(\Graph)$ is an NPC 
	(non-positively curved) space: there exists a 
	metric whose curvature (defined appropriately at the
	non-manifold points) is never positive. 
\item
	The fundamental groups of these spaces are all 
	torsion-free. In other words, if a loop in 
	$\conf^N(\Graph)$ cannot be shrunk to a point, 
	then neither can any multiple of the loop. This 
	property is also true for configuration spaces of 
	$\real^2$, but {\em not} for $\conf^N(S^2)$ --- robots
	on a two-dimensional sphere.  
\item
	The fundamental groups of these spaces have solvable
	word problem, which means that there is an algorithm
	which can be used to decide whether any given loop in 
	$\conf^N(\Graph)$ can be shrunk to a point in $\conf^N(\Graph)$.
\item 
	The fundamental group $\pi_1(\conf^N(\Graph))$ can be
	realized as a {\em graph of groups} --- that is, as a 
	collection of groups (abstractly thought of as vertices
	in some graph) which are pairwise glued together 
	or ``amalgamated'' along a network of subgroups 
	(each gluing represented by an edge between vertices).
\end{enumerate}

For those not familiar with these more subtle features of 
topological spaces and their fundamental groups, we would offer 
the examples presented here
as an excellent concrete manifestation of these properties.
One can easily explain what a configuration space of robots
is to a high-school class\footnote{Both authors have 
done this on several occasions with positive results.},
and, upon demonstrating that, e.g., two robots on 
$K_5$ yields a genus six surface, one has an excellent 
demonstration of the sublime nature of configuration 
spaces and of topology in general.


\begin{thebibliography}{1}

\bibitem{Abr00}
A.~Abrams.
\newblock {\em Configuration Spaces and Braid Groups of Graphs}.
\newblock PhD thesis, UC Berkeley, 2000.

\bibitem{TC98}	
A.~Abrams, H.~Landau, Z.~Landau, J.~Pommersheim, and E.~Zaslow.
\newblock Evasive random walks. 
\newblock Preprint, submitted for publication, 2000.

\bibitem{Bir74}
J.~Birman.
\newblock {\em Braids, Links, and Mapping Class Groups}.
\newblock Princeton University Press, Princeton, N.J., 1974.

\bibitem{BottTaub}
R.~Bott and C.~Taubes.
\newblock On the self-linking of knots: topology and physics.
\newblock {\em J. Math. Phys.} 35(10), 5247--5287, 1994.

\bibitem{Can88}
J.~Canny. 
\newblock {\it The Complexity of Robot Motion Planning}. 
\newblock {MIT Press, Cambridge, MA}, 1988.

\bibitem{agvhandbook}
G. Castleberry.
\newblock {\em The AGV Handbook}.
\newblock Braun-Brumfield, Ann Arbor, MI, 1991.

\bibitem{G:birman}
R.~Ghrist.
\newblock Configuration spaces of graphs in robotics.
\newblock In {\em Braids, Links, and Mapping Class Groups: the Proceedings of 
		Joan Birman's 70th Birthday}, AMS/IP Studies in
		Mathematics vol. 19, 31--41, 2000. 
\newblock ArXiv preprint {\bf math.GT/9905023}. 


\bibitem{GK98}
R.~Ghrist and D.~Koditschek.
\newblock Safe cooperative robot dynamics via dynamics on graphs.
\newblock In Y.~Nakayama, editor, {\em Proceedings of the Eigth International
  Symposium on Robotics Research}, pages 81--92. Springer-Verlag, 1998.

\bibitem{GK00}
R.~Ghrist and D.~Koditschek.
\newblock Safe, cooperative robot dynamics on graphs.
\newblock Preprint, submitted for publication, 2000.
\newblock ArXiv preprint {\bf  cs.RO/0002014}.

\bibitem{Got88}
D. Gottlieb.
\newblock Topology and the robot arm.
\newblock {\em  Acta Appl. Math.}, 11(2): 117--121, 1988.

\bibitem{KM96}
M. Kapovich and J. J. Millson.
\newblock The symplectic geometry of polygons in Euclidean space.
\newblock {\em J. Differential Geom.}, 44(3):479--513, 1996.

\bibitem{KR90}
D.~Koditschek and E.~Rimon.
\newblock Robot navigation functions on manifolds with boundary.
\newblock {\em Adv. in Appl. Math.}, 11(4):412--442, 1990.

\bibitem{latombe.book}
J.-C. Latombe.
\newblock {\em Robot Motion Planning}.
\newblock Kluwer Academic Press, Boston, MA, 1991.

\bibitem{Wal85}
K. Walker.
\newblock {\em Configuration Spaces of Linkages}. 
\newblock {Undergraduate thesis, Princeton University, 1985.}

\bibitem{Wee85}
J. Weeks.
\newblock {\em The Shape of Space}.
\newblock Marcel Dekker Inc., New York, 1985.




\end{thebibliography}


\end{document}